\documentclass{article}
\usepackage{amssymb,amsmath,bm,geometry,graphics,graphicx,url,color,hyperref}
\def\no{\noindent}
\def\pmatrix{\left(\begin{array}}
\def\endpmatrix{\end{array}\right)}
\def\udots{\reflectbox{$\ddots$}}  

\def\RR{\mathbb{R}}

\def\I{{\cal I}}

\def\P{{\cal P}}

\def\dd{\mathrm{d}}

\def\diag{\mathrm{diag}}

\newtheorem{theo}{Theorem}
\newtheorem{lem}{Lemma}

\newtheorem{rem}{Remark}
\newtheorem{defi}{Definition}
\def\proof{\noindent\underline{Proof}\quad}
\def\QED{\mbox{~$\Box{~}$}}

\def\bfc{{\bm{c}}}

\def\bfe{{\bm{e}}}

\def\bfzero{{\bm{0}}}
\def\bfuno{{\bm{1}}}

\def\bfgamma{{\bm{\gamma}}}
\def\bfGamma{{\bm{\Gamma}}}

\def\bfeta{{\bm{\eta}}}

\def\bfrho{{\bm{\rho}}}

\def\hb{{\hat b}}
\def\hc{{\hat c}}

\begin{document}

\title{Arbitrarily high-order energy-preserving methods for simulating the gyrocenter dynamics of charged particles}

\author{Luigi Brugnano$\,^{\rm a}$  \and  Felice Iavernaro$\,^{\rm b}$ \and Ruili Zhang$\,^{\rm c}$}

\date{\small a) Dipartimento di Matematica e Informatica ``U.\,Dini'', Universit\`a di Firenze, 50134 Firenze, Italy.\\ \url{luigi.brugnano@unifi.it}\\
b) Dipartimento di Matematica, Universit\`a di Bari, 70125 Bari, Italy.\\
\url{felice.iavernaro@uniba.it}\\
c) School of Science, Beijing Jiaotong University, Beijing, 100044, China.\\
\url{rlzhang@ustc.edu.cn}}

\maketitle

\begin{abstract}  Gyrocenter dynamics of charged particles plays a fundamental role in plasma physics. In particular, accuracy and conservation of energy are important features for correctly performing long-time simulations. For this purpose, we here propose arbitrarily high-order energy conserving methods for its simulation. The analysis and the efficient implementation of the methods are fully described, and some  numerical tests are reported.

\medskip
\no{\bf Keywords:~ } gyrocenter dynamics,  energy-conserving methods, Poisson problems, line integral methods, Hamiltonian Boundary Value Methods, HBVMs.

\medskip
\no{\bf MSC:~}  65P10, 65L05.

\end{abstract}

\section{Introduction}\label{intro}

The motion of an electrically charged particle, such as an electron or ion, in a plasma  plunged in an electro-magnetic field can be treated as the superposition of a relatively fast circular motion (gyromotion) around a point called the gyrocenter, and a relatively slow drift of this point. Its  numerical simulation is an important issue in plasma physics but brings a number of difficulties.  Not only does the problem exhibit different time-scales due to the fast gyromotion and the slow gyrocenter motion, but it also is multi-scale in space, since the fast gyromotion is small-scale, whereas the slow gyrocenter motion is large-scale. However, when we pay attention to the orbit of a charged particle, the fast gyromotion usually has not much information and, moreover, it is time consuming in numerical simulations. Consequently, we shall focus on the orbits of the gyrocenter, thus only considering the gyrocenter dynamics of a charged particle.

The mathematical description of the gyrocenter dynamics of a charged particle has been at first given by R.G.\,Littlejohn \cite{LJ1983}. In more details, we shall consider the case of time-independent magnetic and electric fields, respectively given by
\begin{equation}\label{BE}
B(x)=\nabla\times A(x)\qquad \mbox{and} \qquad E(x)=-\nabla \phi(x), \qquad x\in\RR^3,
\end{equation}
with $\phi(x)$ the electric potential and
\begin{equation}\label{A}
A(x)\equiv\pmatrix{c}A_1(x)\\ A_2(x)\\ A_3(x)\endpmatrix
\end{equation}
the magnetic vector potential. In such a case, by setting $x=(x_1,x_2,x_3)^\top\in\RR^3$ the position of the gyrocenter, defining\,\footnote{Hereafter, $\|\cdot\|$ will denote the Euclidean norm.}
\begin{equation}\label{b}
b(x)\equiv \pmatrix{c}b_1(x)\\b_2(x)\\b_3(x)\endpmatrix :=\dfrac{B(x)}{\|B(x)\|}
\end{equation}
as the unit vector in the direction of the magnetic field, along which the gyrocenter moves with velocity $u$, and setting
\begin{equation}\label{y}
y:=\pmatrix{c}x\\ u\endpmatrix \in\RR^4,
\end{equation}
the equations describing the gyrocenter dynamics of the charged particle can be written as
\begin{equation}\label{ghs}
\dot y =K(y)^{-1}\nabla H(y),
\end{equation}
where:
{ 
\begin{itemize}

\item the {\em Hamiltonian function} $H(y)$ is given by
\begin{equation}\label{H}
H(y)=\dfrac{1}{2}u^2+\mu  \|B(x)\|+\phi(x),
\end{equation}
with $\mu$ a given constant representing an adiabatic invariant;

\item $K^{-1}(y)$ is the following skew-symmetric matrix,
\begin{equation}\label{K1}
K^{-1}(y)=\frac{1}{|b(x)^\top a(y)|}\pmatrix{cccc}
0 & -b_{3}(x) & b_{2}(x) & a_1(y)\\
b_{3}(x) & 0 & -b_{1}(x) &  a_2(y)\\
-b_{2}(x) & b_{1}(x) & 0 & a_3(y)\\
-a_1(y) & -a_2(y) & -a_3(y) & 0
\endpmatrix,
\end{equation}
with $b(x)$ the vector defined in (\ref{b}) and
\begin{equation}\label{aij}
a(y) \equiv \pmatrix{c} a_1(y) \\ a_2(y)\\ a_3(y)\endpmatrix  := \nabla\times \left[ A(x)+u\, b(x)\right]  = B(x) + u\nabla\times b(x).
\end{equation}
\end{itemize}
}

Consequently, {  without loss of generality,} the problem we want to solve can be assumed in the form  
\begin{equation}\label{poisson}
\dot y = S(y) \nabla H(y) =: f(y), \qquad t\in[0,h],\qquad y(0)=y_0, \qquad S(y)=-S(y)^\top,
\end{equation}
where $h$ is the timestep used {  and, for sake of brevity, the notation $S(y)=K^{-1}(y)$ has been used}. 
Clearly, for such a problem, the energy $H$ is conserved along the solution, since
\begin{equation}\label{dH}
\dot H(y) = \nabla H(y)^\top\dot y = \nabla H(y)^\top S(y) \nabla H(y) = 0,
\end{equation}
due to the skew-symmetry of $S(y)$. Energy-conserving methods for problems in this form have been already considered, e.g., in \cite{CH2011,BCMR2012,ZLQT2019}. Here,  the derivation and analysis of the methods will be done within the framework of {\em line integral methods}, namely methods defined by a suitable smooth path $\sigma$ such that:
\begin{equation}\label{lim12}
\sigma(0) = y_0, \qquad y_1 :=\sigma(h),
\end{equation}
with $y_1$ the approximation to $y(h)$, and\,\footnote{The denomination {\em line integral methods} stems from the fact that one requires the vanishing of the line integral along the path $\sigma$, as prescribed by (\ref{lim3}).}
\begin{equation}\label{lim3}
H(y_1) -H(y_0) = H(\sigma(h))-H(\sigma(0)) = h\int_0^1 \nabla H(\sigma(ch))^\top\dot\sigma(ch)\dd c = 0.
\end{equation}
In such a case, the conservation of $H$ will no more require the integrand to be identically zero, as is the case for $y(t)$ (see (\ref{dH})), thus leaving much more freedom in the choice of the path $\sigma$. Line integral methods have been at first introduced for deriving energy-conserving methods for Hamiltonian problems  \cite{IP2007,IP2008,IT2009,BIS2009,BIT2009,BIT2010,BIT2012}, resulting in the class of energy-conserving methods named {\em Hamiltonian Boundary Value Methods (HBVMs)}, and have then been developed along several directions (see, e.g.,  the monograph \cite{LIMbook2016} or the review papers \cite{BI2018,BFCI2019}), including their use as spectral methods in time \cite{BMR2018,BIMR2019,ABI2019,ABI2019_1}.

With these premises, in Section~\ref{cont} we derive and analyse,   at a continuous level, polynomial approximations able to exactly conserve the energy $H$;  fully discrete methods are then derived and studied in Section~\ref{discre}, where we also sketch the efficient solution of the generated discrete problems; finally,  a few numerical tests are reported in Section~\ref{num}, along with some concluding remarks.

\section{Derivation of the method}\label{cont}

The approach used to derive the method is based on a suitable modification of the arguments used in \cite{BIT2012,BMR2019}. To begin with, let us expand the gradient of $H$ along the orthonormal Legendre polynomial basis on $[0,1]$:
\begin{equation}\label{leg}
P_i\in\Pi_i, \qquad \int_0^1 P_i(x)P_j(x)\dd x =\delta_{ij}. \qquad i,j=0,1,\dots,
\end{equation}
where, as usual, $\Pi_i$ is the set of polynomials of degree at most $i$, and $\delta_{ij}$ is the Kronecker delta. Moreover, for sake of simplicity, we shall hereafter assume that both $\nabla H(y(t))$ and $S(y(t))$ admit a Taylor expansion at $t=0$. Consequently, one has
\begin{equation}\label{gammaj}
\nabla H(y(ch)) = \sum_{j\ge0} P_j(c)\gamma_j(y), \qquad \gamma_j(y) = \int_0^1P_j(\tau)\nabla H(y(\tau h))\dd\tau, \quad j=0,1,\dots,
\end{equation}
and, therefore, the equation in (\ref{poisson}) can be rewritten as
\begin{equation}\label{poisson1}
\dot y(ch) = S(y(ch))\sum_{j\ge0} P_j(c)\gamma_j(y), \qquad c\in[0,1].
\end{equation}
Next, let us consider the expansions
\begin{equation}\label{roij}
P_j(c)S(y(ch)) = \sum_{i\ge0} P_i(c)\rho_{ij}(y), \quad \rho_{ij}(y) = \int_0^1P_i(\tau)P_j(\tau)S(y(\tau h))\dd\tau \equiv \rho_{ji}(y), \quad i,j=0,1,\dots.
\end{equation}
Eventually, we obtain that the equation in (\ref{poisson}) can be rewritten, by virtue of (\ref{gammaj})--(\ref{roij}), as
\begin{equation}\label{Gammai}
\dot y(ch) = \sum_{i,j\ge0} P_i(c) \rho_{ij}(y)\gamma_j(y) \equiv \sum_{i\ge0} P_i(c)\Gamma_i(y), \qquad c\in[0,1],
\qquad \Gamma_i(y) = \sum_{j\ge0} \rho_{ij}(y)\gamma_j(y),
\end{equation}
from which, integrating term by term, and imposing the initial condition $y(0)=y_0$, one has:
\begin{equation}\label{yt}
y(ch) = y_0 + h\sum_{i\ge0}\int_0^c P_i(x)\dd x \, \Gamma_i(y), \qquad c\in[0,1].
\end{equation}
\begin{rem}\label{fouGi}
It is clear, from (\ref{Gammai}), that $\Gamma_i(y)$ is nothing but the $i$-th Fourier coefficient in the expansion of the right-hand side of (\ref{poisson})  along the Legendre basis (\ref{leg}), i.e.,
\begin{equation}\label{fyexp}
f(y(ch)) = \sum_{i\ge0} P_i(c)\Gamma_i(y), \qquad c\in[0,1].
\end{equation}
\end{rem}

We need the following preliminary results.

\begin{lem}\label{lem1} Let $g:[0,h]\rightarrow V$, with $V$ a vector space, admit a Taylor expansion at 0. Then,
$$\int_0^1 P_j(c)g(ch)\dd c = O(h^j), \qquad j=0,1,\dots.$$
\end{lem}
\proof See \cite[Lemma~1]{BIT2012}.\,\QED\medskip

\begin{lem}\label{lem2} $\gamma_j(y)=O(h^j)$, $\rho_{ij}(y)=-\rho_{ij}(y)^\top = O(h^{|j-i|})$.\end{lem}
\proof The first point follows from Lemma~\ref{lem1}, the second from \cite[Corollary\,1]{BMR2019} and the skew-symmetry of $S(y)$.\,\QED\medskip

\begin{lem}\label{lem3}$\Gamma_i(y) = O(h^i)$.\end{lem}
\proof The statement follows from Lemma~\ref{lem1}, by taking into account  (\ref{fyexp}).\,\QED\bigskip

In order to obtain a polynomial approximation of degree $s$ to $y$, we truncate the series in (\ref{gammaj}) and (\ref{roij}) to finite sums, thus getting
\begin{equation}\label{sig1}
\dot\sigma(ch) = \sum_{i=0}^{s-1} P_i(c) \Gamma_i^s(\sigma), \qquad c\in[0,1], \qquad \Gamma_i^s(\sigma) = \sum_{j=0}^{s-1}\rho_{ij}(\sigma)\gamma_j(\sigma),\quad i=0,\dots,s-1,
\end{equation}with $\gamma_j(\sigma)$ and $\rho_{ij}(\sigma)$ formally still defined by (\ref{gammaj}) and (\ref{roij}), respectively, upon  replacing $y$ by $\sigma$. It is straightforward to verify that Lemma~\ref{lem2} continues to hold for the Fourier coefficients $\rho_{ij}(\sigma)$ and $\gamma_j(\sigma)$. Moreover, similarly as for Lemma~\ref{lem3}, now one has the following result which holds true for $\Gamma_i^s(\sigma)$.

\begin{lem}\label{lem4}For all\, $ i=0,\dots,s-1$, one has: $\Gamma_i^s(\sigma) = O(h^i), ~ \Gamma_i(\sigma)-\Gamma_i^s(\sigma) = O(h^{2s-i})$.\end{lem}
\proof One has, by virtue of Lemma~\ref{lem2}, and for $i=0,\dots,s-1$:
$$\Gamma_i^s(\sigma) = \sum_{j=0}^{s-1} \rho_{ij}(\sigma)\gamma_j(\sigma) = \sum_{j=0}^i\underbrace{O(h^{i-j})O(h^j)}_{=O(h^i)}+
\sum_{j=i+1}^{s-1}\underbrace{O(h^{j-i})O(h^j)}_{=O(h^{2j-i})} = O(h^i).$$
Similarly, for $i<s$, one obtains:
$$\Gamma_i(\sigma)-\Gamma_i^s(\sigma)=\sum_{j\ge s} \rho_{ij}(\sigma)\gamma_j(\sigma) = \sum_{j\ge s}\underbrace{O(h^{j-i})O(h^j)}_{=O(h^{2j-i})} = O(h^{2s-i}).\,\QED$$

Integrating term by term (\ref{sig1}), and imposing the initial condition $\sigma(0)=y_0$, one then obtains the polynomial approximation of degree $s$
\begin{equation}\label{sig}
\sigma(ch) = y_0 + h\sum_{i=0}^{s-1} \int_0^c P_i(x)\dd x\, \Gamma_i^s(\sigma), \qquad c\in[0,1],
\end{equation}
with the approximation to $y(h)$ given, by considering that $\int_0^1P_i(x)\dd x=\delta_{i0}$ (see (\ref{leg})), by:
\begin{equation}\label{y1}
y_1:=\sigma(h) = y_0+h\Gamma_0^s(\sigma).
\end{equation}
In so doing, the requirements (\ref{lim12}) are fulfilled.

\subsection{Analysis}\label{anal}

Let us now study the properties of the approximation procedure (\ref{sig1})--(\ref{y1}). To begin with, let us now prove that also the requirement (\ref{lim3}) holds true, besides (\ref{lim12}).

\begin{theo}\label{cons1} For the  polynomial approximation defined by (\ref{sig1})--(\ref{y1}), one has: $H(y_1)=H(y_0)$,  i.e., the method is energy-conserving. \end{theo}
\proof In fact, according to (\ref{lim3}) one obtains:
\begin{eqnarray*}
\lefteqn{
H(y_1)-H(y_0)~=~H(\sigma(h))-H(\sigma(0)) ~=~ h\int_0^1 \nabla H(\sigma(ch))^\top\dot\sigma(ch)\dd c} \\
&=& h\int_0^1 \nabla H(\sigma(ch))^\top\sum_{i=0}^{s-1} P_i(c) \Gamma_i^s(\sigma) \dd c~=~ h \sum_{i=0}^{s-1}\int_0^1 \nabla H(\sigma(ch))^\top P_i(c) \dd c\,\Gamma_i^s(\sigma)\\
&=& h\sum_{i=0}^{s-1} \gamma_i(\sigma)^\top \Gamma_i^s(\sigma) ~=~ h\sum_{i,j=0}^{s-1}  \gamma_i(\sigma)^\top \rho_{ij}(\sigma)\gamma_j(\sigma) ~=~ 0,
\end{eqnarray*}
due to the skew-symmetry of $\rho_{ij}(\sigma).\,\QED$\bigskip

Next, let us discuss the accuracy of the approximation (\ref{y1}) to $y(h)$, adapting the arguments of \cite[Theorem\,1]{BIT2012}. Preliminarily, let us denote by $y(t,t^*,y^*)$ the solution of the problem (see (\ref{poisson}))
\begin{equation}\label{fy}
\dot y = f(y), \qquad t\in[t^*,h], \qquad y(t^*)=y^*,
\end{equation}
also recalling the following well-known perturbation results,
\begin{equation}\label{perturb}
\frac{\partial}{\partial y^*} y(t,t^*,y^*) = \Phi(t,t^*), \qquad  \frac{\partial}{\partial t^*} y(t,t^*,y^*) = -\Phi(t,t^*)f(y^*),
\end{equation}
where $\Phi(t,t^*)$ is the fundamental matrix solution of the variational problem associated to (\ref{fy}).

\begin{theo}\label{y1ord}  For the polynomial approximation defined by (\ref{sig1})--(\ref{y1}), one has: $y_1-y(h)=O(h^{2s+1})$, where $y(t)$ is the solution of problem (\ref{poisson}).\footnote{I.e., the method has order $2s$.} \end{theo}
\proof By virtue of Lemmas~\ref{lem1}, \ref{lem3}, and \ref{lem4}, from  (\ref{fy})-(\ref{perturb}), and considering that, according to (\ref{fyexp}),
$$f(\sigma(ch))=\sum_{j\ge0}P_j(c)\Gamma_i(\sigma),\qquad c\in[0,1],$$ one has:
\begin{eqnarray*}
\lefteqn{
y_1-y(h)~=~y(h,h,y_1)-y(h,0,y_0) ~=~ y(h,h,\sigma(h))-y(h,0,\sigma(0)) ~=~ \int_0^h \frac{\dd}{\dd t} y(h,t,\sigma(t))\dd t}\\
             &=& \int_0^h \left[\left.\frac{\partial}{\partial t^*} y(h,t^*,\sigma(t))\right|_{t^*=t}+ \left.\frac{\partial}{\partial y^*} y(h,t,y^*)\right|_{y^*=\sigma(t)}\dot\sigma(t)\right]\dd t\\
             &=& \int_0^h \left[-\Phi(h,t)f(\sigma(t))+\Phi(h,t)\dot\sigma(t)\right]\dd t ~=~ h\int_0^1 \Phi(h,ch)\left[\dot\sigma(ch)-f(\sigma(ch))\right]\dd c\\
             &=& h\sum_{i=0}^{s-1} \underbrace{\int_0^1 P_i(c)\Phi(h,ch)\dd c}_{=O(h^i)}\underbrace{\left[ \Gamma_i^s(\sigma)-\Gamma_i(\sigma)\right]}_{=O(h^{2s-i})} - h\sum_{i\ge s}\underbrace{\int_0^1 P_i(c)\Phi(h,ch)\dd c}_{=O(h^i)}\underbrace{\Gamma_i(\sigma)}_{O(h^i)}
             ~=~O(h^{2s+1}).\,\QED
\end{eqnarray*}

\section{Discretization}\label{discre}
As is clear, the polynomial approximation (\ref{sig1})--(\ref{y1}) does not yet provide a {\em numerical method}. As matter of fact, quoting \cite[p.\,521]{DB2008}, {\em ``as is well known, even many relatively simple integrals cannot be expressed in finite terms of elementary functions, and thus must be evaluated by numerical methods''}. In the present setting, this means that we need to approximate the Fourier coefficients
\begin{equation}\label{rog}
\rho_{ij}(\sigma)=\int_0^1P_i(c)P_j(c)S(\sigma(ch))\dd c,\quad \gamma_j(\sigma)=\int_0^1P_j(c)\nabla H(\sigma(ch))\dd c, \quad i,j=0,\dots,s-1,
\end{equation} by using suitable quadratures.
For this purpose, we shall use:
\begin{itemize}
\item a Gauss-Legendre formula of order $2k_1$, with $k_1\ge s$, to approximate the former coefficients, whose abscissae and weights shall be  denoted by
\begin{equation}\label{hbc}
\hc_\ell,\,\hb_\ell, \qquad \ell=1,\dots,k_1,
\end{equation}respectively;
\item a Gauss-Legendre formula of order $2k_2$, with $k_2\ge s$, to approximate the latter coefficients, whose abscissae and weights shall be   denoted by
\begin{equation}\label{bc}
c_\ell,\,b_\ell, \qquad \ell=1,\dots,k_2,
\end{equation}respectively.
\end{itemize}
We recall that, due to the symmetry of the abscissae (\ref{hbc}) and (\ref{bc}), one has:
\begin{equation}\label{simm1}
\hat c_\ell = 1-\hat c_{k_1-\ell+1}, \quad \hat b_\ell = \hat b_{k_1-\ell+1}, \quad \ell=1,\dots, k_1,
\end{equation}
and
\begin{equation}\label{simm2}
c_\ell = 1-c_{k_2-\ell+1}, \quad b_\ell = b_{k_2-\ell+1}, \qquad \ell=1,\dots, k_2.
\end{equation}
As a result, in place of the polynomial $\sigma$ defined at (\ref{sig1})--(\ref{sig}), we shall have a, generally different, polynomial $u\in\Pi_s$, such that:
\begin{equation}\label{u1}
\dot u(ch) = \sum_{i=0}^{s-1} P_i(c) \hat\Gamma_i, \qquad c\in[0,1], \qquad \hat\Gamma_i = \sum_{j=0}^{s-1}\hat\rho_{ij}\hat\gamma_j,\quad i=0,\dots,s-1,\end{equation}
with the approximate Fourier coefficients
\begin{equation}\label{hrogam}
\hat\rho_{ij} = \sum_{\ell=1}^{k_1} \hb_\ell P_i(\hc_\ell)P_j(\hc_\ell)S(u(c_\ell h))\equiv \hat\rho_{ji}, \quad
\hat\gamma_j = \sum_{\ell=1}^{k_2} b_\ell P_j(c_\ell)\nabla H(u(c_\ell h)), \quad i,j=0,\dots,s-1,
\end{equation}
in place of (\ref{rog}). Consequently, one obtains
\begin{equation}\label{u}
u(ch) = y_0+h\sum_{i=0}^{s-1} \int_0^c P_i(x)\dd x\, \hat\Gamma_i, \qquad c\in[0,1],
\end{equation}
with the new approximation given by (see (\ref{y1}))
\begin{equation}\label{hy1}
y_1 := u(h) \equiv y_0+h\hat\Gamma_0.
\end{equation}
For later use, let us recall that the errors in the quadratures (\ref{hrogam}) are given by (see (\ref{rog})),\footnote{Hereafter, $S\in\Pi_\nu$ means that the entries of matrix $S$ in (\ref{poisson}) are polynomials of degree at most $\nu$ in the argument.}
\begin{equation}\label{hdeltaij}
\hat\rho_{ij}-\rho_{ij}(u) = \hat\Delta_{ij}(h) \equiv
\left\{\begin{array}{cl}0, &\mbox{if}\quad S\in\Pi_\nu \mbox{~~with~~} \nu< [2(k_1-s)+1]/s,\\[2mm] O(h^{2k_1-i-j}), &\mbox{otherwise.}\end{array}\right.
\end{equation}
and 
\begin{equation}\label{deltaj}
\hat\gamma_j-\gamma_j(u) = \Delta_j(h) \equiv \left\{\begin{array}{cl}0, &\mbox{if}\quad H\in\Pi_\nu \mbox{~~with~~} \nu\le 2k_2/s,\\[2mm] O(h^{2k_2-j}), &\mbox{otherwise.}\end{array}\right.
\end{equation}

\begin{rem}\label{k1k2} As is clear from (\ref{hdeltaij}) and (\ref{deltaj}), when both the entries of the matrix $S$ and the energy $H$ in (\ref{poisson}) are polynomials, we can exactly compute the Fourier coefficients, by choosing $k_1$ and $k_2$ large enough, so that $u\equiv\sigma$. This, in turn, will be not a big computational issue since, as we shall see later,  the discrete problem will always have dimension $s$, {\em independently} of $k_1$ and $k_2$.\end{rem}

It is quite straightforward to prove the following properties, representing the discrete counterpart of Lemma~\ref{lem2}:
\begin{equation} \label{discpro}
\forall k_1,k_2\ge s: \quad \hat\gamma_j = O(h^j), \quad \hat\rho_{ij}=-\hat\rho_{ij}^\top=O(h^{|j-i|}), \quad i,j=0,\dots,s-1.
\end{equation}
Moreover, the following result holds also true.\footnote{It represents the discrete counterpart of Lemma~\ref{lem4}.}

\begin{lem}\label{lem5}
$\forall k_1,k_2\ge s:\quad \hat\Gamma_i = O(h^i), \quad \hat\Gamma_i-\Gamma_i(u) = O(h^{2s-i}), \quad i=0,\dots,s-1.$
\end{lem}
\proof The first part of the statement follows from the second one and from Lemma~\ref{lem3}, by considering that $i<s$. Next, by virtue of Lemma~\ref{lem4}, one has:
$$ \hat\Gamma_i-\Gamma_i(u) =  \hat\Gamma_i- \Gamma_i^s(u)+\Gamma_i^s(u)-\Gamma_i(u) =  \hat\Gamma_i- \Gamma_i^s(u)+O(h^{2s-i}).$$
Moreover, from (\ref{sig1}) and (\ref{hdeltaij})-(\ref{deltaj}), one has:
\begin{eqnarray*}
\hat\Gamma_i- \Gamma_i^s(u) &=& \sum_{j=0}^{s-1} \hat\rho_{ij}\hat\gamma_j-\rho_{ij}(u)\gamma_j(u)
= \sum_{j=0}^{s-1} (\rho_{ij}(u)+\hat\Delta_{ij}(h))(\gamma_j(u)+\Delta_j(u))-\rho_{ij}(u)\gamma_j(u)\\
&=&\sum_{j=0}^{s-1} \rho_{ij}(u)\Delta_j(u)+\hat\Delta_{ij}(h)\gamma_j(u)+\hat\Delta_{ij}(h)\Delta_j(h).
\end{eqnarray*}
Considering that:
\begin{itemize}
\item $\rho_{ij}(u)\Delta_j(u) = O(h^{|j-i|+2k_2-j})$ ~and~ $|j-i|+2k_2-j\ge 2k_2-i\ge 2s-i$;
\item $\hat\Delta_{ij}(h)\gamma_j(u) = O(h^{2k_1-i-j+j})$ ~and~ $2k_1-i\ge 2s-i$;
\item $\hat\Delta_{ij}(h)\Delta_j(h) = O(h^{2k_1-i-j+2k_2-j})$ ~and~ $2(k_1+k_2-j)-i\ge 2s-i$;
\end{itemize}
the statement then follows.\,\QED\smallskip

\begin{defi}\label{lim} We call {\em discrete line integral method with parameters $k_1,k_2,s$}, in short\, {\em LIM$(k_1,k_2,s)$}, the numerical method defined by  (\ref{hbc})--(\ref{hy1}).
\end{defi}

\subsection{Analysis}\label{anal2}

Let us now study the properties of the  LIM$(k_1,k_2,s)$ method (\ref{hbc})--(\ref{hy1}). We start discussing the conservation of energy.

\begin{theo}\label{cons2} For the LIM$(k_1,k_2,s)$ method, {  with $k_1,k_2\ge s$,} one has:
$$H(y_1)-H(y_0) = \left\{ \begin{array}{cl} 0, &\mbox{~if\quad} H\in\Pi_\nu \mbox{~~with~~} \nu\le2k_2/s,\\[2mm]
O(h^{2k_2+1}), &\mbox{otherwise.}\end{array}\right.$$\end{theo}
\proof In fact, by taking into account (\ref{gammaj}) and (\ref{u1})--(\ref{discpro}), one obtains:
\begin{eqnarray*}
\lefteqn{
H(y_1)-H(y_0)~=~H(u(h))-H(u(0)) ~=~ h\int_0^1 \nabla H(u(ch))^\top\dot u(ch)\dd c} \\
&=& h\int_0^1 \nabla H(u(ch))^\top\sum_{i=0}^{s-1} P_i(c) \hat\Gamma_i \dd c~=~ h \sum_{i=0}^{s-1}\int_0^1 \nabla H(u(ch))^\top P_i(c) \dd c\,\hat\Gamma_i\\
&=& h\sum_{i=0}^{s-1} \gamma_i(u)^\top \hat\Gamma_i ~=~ h\sum_{i,j=0}^{s-1}  \gamma_i(u)^\top \hat\rho_{ij}\hat\gamma_j
~=~ h\sum_{i,j=0}^{s-1}  \gamma_i(u)^\top \hat\rho_{ij}\left(\gamma_j(u)+\Delta_j(h)\right)\\
&=& h\sum_{i,j=0}^{s-1}  \gamma_i(u)^\top \hat\rho_{ij}\Delta_j(h).
\end{eqnarray*}
When $H\in\Pi_\nu$, with $\nu\le 2k_2/s$, the quadrature error $\Delta_j(h)=0$. Conversely, one has:
$${   h\sum_{i,j=0}^{s-1}  \gamma_i(u)^\top \hat\rho_{ij}\Delta_j(h)} = h\sum_{i,j=0}^{s-1}  \underbrace{\gamma_i(u)^\top}_{=O(h^i)} \overbrace{\hat\rho_{ij}}^{=O(h^{|j-i|})}\underbrace{\Delta_j(h)}_{=O(h^{2k_2-j})} = O(h^{2k_2+1}).\,\QED$$


Concerning the order of accuracy, the following result holds true.

 \begin{theo}\label{hy1ord}  For the LIM$(k_1,k_2,s)$ method, if $k_1,k_2\ge s$ one has: $y_1-y(h)=O(h^{2s+1})$, where $y(t)$ is the solution of problem (\ref{poisson}).\end{theo}
\proof Following similar steps as those used in the proof of Theorem~\ref{y1ord}, and taking into account (\ref{u1})--(\ref{discpro}) and Lemma~\ref{lem5}, one has:
\begin{eqnarray*}
\lefteqn{y_1-y(h)~=~y(h,h,y_1)-y(h,0,y_0) ~=~ y(h,h,u(h))-y(h,0,u(0)) ~=~ \int_0^h \frac{\dd}{\dd t} y(h,t,u(t))\dd t}\\
             &=& \int_0^h \left[\left.\frac{\partial}{\partial t^*} y(h,t^*,u(t))\right|_{t^*=t}+ \left.\frac{\partial}{\partial y^*} y(h,t,y^*)\right|_{y^*=u(t)}\dot u(t)\right]\dd t\\
             &=& \int_0^h \left[-\Phi(h,t)f(u(t))+\Phi(h,t)\dot u(t)\right]\dd t ~=~ h\int_0^1 \Phi(h,ch)\left[\dot u(ch)-f(u(ch))\right]\dd c\\
             &=& h\sum_{i=0}^{s-1} \underbrace{\int_0^1 P_i(c)\Phi(h,ch)\dd c}_{=O(h^i)}\underbrace{\left[ \hat\Gamma_i-\Gamma_i(u)\right]}_{=O(h^{2s-i})} - h\sum_{i\ge s}\underbrace{\int_0^1 P_i(c)\Phi(h,ch)\dd c}_{=O(h^i)}\underbrace{\Gamma_i(u)}_{O(h^i)}
             ~=~O(h^{2s+1}).\,\QED
\end{eqnarray*}\medskip

{ 
\begin{rem}\label{k2} From the result of the last two theorems, one deduces that a LIM$(k_1,k_2,s)$ method:
\begin{itemize}
\item has order $2s$, provided that $k_1,k_2\ge s$;
\item is energy-conserving, when the Hamiltonian $H$ in (\ref{poisson}) is a polynomial of degree not larger that $2k_2/s$.
\end{itemize}
We observe that, even when $H$ is not a polynomial, a {\em practical} energy conservation can always be gained by choosing $k_2$ large enough, so that the $O(h^{2k_2+1})$ energy error falls below the round-off error level: for this reason, in the sequel we shall make no distinction between the two cases.
\end{rem}
}

Let us now derive a compact formulation of the discrete problem generated by the method. As is clear, it is enough to compute the coefficients $\hat\Gamma_i$, $i=0,\dots,s-1$, of the polynomial approximation (\ref{u}). Let us then define the following vectors and matrices:
\begin{equation}\label{GRgam}
\bfGamma = \pmatrix{c} \hat\Gamma_0\\ \vdots \\ \hat\Gamma_{s-1}\endpmatrix, \qquad
\bfrho = \pmatrix{ccc} \hat\rho_{00} & \dots &\hat\rho_{0,s-1}\\
\vdots & & \vdots\\
\hat\rho_{s-1,0} & \dots & \hat\rho_{s-1,s-1}\endpmatrix,\qquad
\bfgamma = \pmatrix{c}\hat\gamma_0\\ \vdots \\ \hat\gamma_{s-1}\endpmatrix,
\end{equation}

\begin{equation}\label{uhu}
\hat\bfc = \pmatrix{c} \hat c_1\\ \vdots \\ \hat c_{k_1}\endpmatrix, \quad
\bfc = \pmatrix{c} c_1\\ \vdots \\ c_{k_2}\endpmatrix, \quad
u(\hat\bfc h) = \pmatrix{c} u(\hat c_1h)\\ \vdots\\ u(\hat c_{k_1}h)\endpmatrix, \quad
u(\bfc h) = \pmatrix{c} u(c_1h)\\ \vdots\\ u(c_{k_2}h)\endpmatrix,
\end{equation}

\begin{equation}\label{SDH}
S(u(\hat\bfc h)) = \pmatrix{ccc} S(u(\hat c_1h))\\ &\ddots\\ &&S(u(\hat c_{k_1}h))\endpmatrix, \qquad
\nabla H(u(\bfc h)) = \pmatrix{c} \nabla H(u(c_1h))\\ \vdots \\ \nabla H(u(c_{k_2}h))\endpmatrix,
\end{equation}

\begin{equation}\label{hPIOm}
\hat\P_s = \pmatrix{c} P_{j-1}(\hat c_i)\endpmatrix, ~\hat\I_s =  \pmatrix{c} \int_0^{\hat c_i}P_{j-1}(x)\dd x\endpmatrix\in\RR^{k_1\times s},\qquad
\hat\Omega = \diag\left(\hat b_1,\,\dots,\,\hat b_{k_1}\right),
\end{equation}

\begin{equation}\label{PIOm}
\P_s = \pmatrix{c} P_{j-1}(c_i)\endpmatrix, ~\I_s =  \pmatrix{c} \int_0^{c_i}P_{j-1}(x)\dd x\endpmatrix\in\RR^{k_2\times s},\qquad
\Omega = \diag\left(b_1,\,\dots,\,b_{k_2}\right),
\end{equation}
and
\begin{equation}\label{uhu}
\hat\bfuno = \pmatrix{c} 1\\ \vdots \\1\endpmatrix\in\RR^{k_1}, \qquad \bfuno = \pmatrix{c} 1\\ \vdots \\1\endpmatrix\in\RR^{k_2}.
\end{equation}
Moreover, hereafter we set $I$ the identity matrix having the same dimension as that of the state vector.
One has, then:
\begin{eqnarray*}
u(\hat\bfc h) &=&\hat\bfuno\otimes y_0 + h\hat\I_s\otimes I\,\bfGamma,\\[1mm]
u(\bfc h)      &=&\bfuno\otimes y_0 + h\I_s\otimes I\,\bfGamma,\\[1mm]
\bfgamma   &=& [\P_s^\top\Omega\otimes I]\,\nabla H(\bfuno\otimes y_0 + h\I_s\otimes I\,\bfGamma),\\[1mm]
\bfrho         &=& [\hat\P_s^\top\hat\Omega\otimes I]\, S( \hat\bfuno\otimes y_0 + h\hat\I_s\otimes I\,\bfGamma)\, [\hat\P_s\otimes I],\\[1mm]
\bfGamma &=& \bfrho\bfgamma.
\end{eqnarray*}
Consequently, one obtains the following discrete problem involving only $\bfGamma$:
\begin{equation}\label{bfGampro}
G(\bfGamma) :=  \bfGamma-[\hat\P_s^\top\hat\Omega\otimes I]\, S( \hat\bfuno\otimes y_0 + h\hat\I_s\otimes I\,\bfGamma)\, [\hat\P_s \P_s^\top\Omega\otimes I] \,\nabla H(\bfuno\otimes y_0 + h\I_s\otimes I\,\bfGamma) = \bfzero.
\end{equation}\medskip

\begin{rem} We observe that, in the case $k_1=s$, we obtain the methods defined in \cite{BCMR2012} (compare (\ref{bfGampro}) with \cite[Eq.\,(36)]{BCMR2012}) which, in turn, are akin to those in \cite{CH2011}. Their derivation, however, is now done within a different framework, yielding more general methods (i.e., using $k_1>s$), though having the same order of accuracy.   Moreover, the line integral approach used to derive $(\ref{bfGampro})$ allows us to easily extend a technique devised for the efficient implementation of line integral methods applied to Hamiltonian problems (see e.g. \cite{BIT2011}). This aspect will be faced in the next section.   
\end{rem}

Another important property of the LIM$(k_1,k_2,s)$ method is symmetry, namely if we apply the method to the equation in (\ref{poisson}),  starting from $y_1$ as defined in (\ref{hy1}) and using a stepsize $-h$, this brings us back to $y_0$. To prove this property, we need the following preliminary result.

\begin{lem}\label{lemsym}
With reference to (\ref{hPIOm})--(\ref{uhu}), let us define the following matrices:
\begin{equation}\label{D}
I_\ell = \pmatrix{ccc}1\\ &\ddots\\ &&1\endpmatrix,~ D_\ell = \pmatrix{ccc}(-1)^0\\ &\ddots\\ &&(-1)^{\ell-1}\endpmatrix, ~ P_\ell = \pmatrix{ccc} &&1\\&\udots\\1\endpmatrix\in\RR^{\ell\times\ell},
\end{equation} and set $\bfe_1\in\RR^s$ the first unit vector. Then,
\begin{eqnarray*}
&&D_\ell^2=P_\ell^2=I_\ell, \qquad P_{k_1}\hat\I_s D_s = \hat\bfuno\bfe_1^\top -\hat\I_s,\qquad P_{k_2}\I_s D_s = \bfuno\bfe_1^\top -\I_s,\\[1mm]
&&P_{k_1}\hat\P_s D_s = \hat\P_s, \quad\, P_{k_2}\P_s D_s = \P_s,\quad\,
P_{k_1}\hat\Omega P_{k_1}=\hat\Omega,\quad\, P_{k_2}\Omega P_{k_2}=\Omega.
\end{eqnarray*}
\end{lem}
\proof See \cite[Lemma~3]{BMR2019}. In particular all the properties, but the first one, derive from the symmetry of the abscissae (and, then, of the weights) (\ref{simm1})-(\ref{simm2}).\,\QED\bigskip

We are now in the position of proving the symmetry of the methods.

\begin{theo}\label{thsym} For all $k_1,k_2\ge s$, the method LIM$(k_1,k_2,s)$ is symmetric.\end{theo}
\proof By using the method on the equation (\ref{poisson}) starting at $y_1$, as defined in (\ref{hy1}), with timestep $-h$ we obtain (see (\ref{bfGampro}))
\begin{equation}\label{bgam}
\bar\bfGamma \equiv \pmatrix{c}
\bar\Gamma_0\\ \vdots \\ \bar\Gamma_{s-1}\endpmatrix
=[\hat\P_s^\top\hat\Omega\otimes I]\, S( \hat\bfuno\otimes y_1 - h\hat\I_s\otimes I\,\bar\bfGamma)\, [\hat\P_s \P_s^\top\Omega\otimes I] \,\nabla H(\bfuno\otimes y_1 - h\I_s\otimes I\,\bar\bfGamma),\end{equation}
with the new approximation given by
~$\bar y_0 = y_1-h\bar\Gamma_0.$~ We have then to show that $\bar y_0 = y_0$. This will follow from the fact that (see (\ref{D}))
$$\bar\Gamma_j = (-1)^j\hat\Gamma_j =: \Gamma_j^*, \quad j=0,\dots,s-1, \qquad \Longleftrightarrow \qquad\bar\bfGamma = D_s\otimes I\, \bfGamma  =:\pmatrix{c} \Gamma_0^*\\ \vdots \\ \Gamma_{s-1}^*\endpmatrix\equiv \bfGamma^*.$$
We shall prove this statement by showing that $\bfGamma^*$ satisfies the same equation (\ref{bgam}) which implicitly defines $\bar\bfGamma$. One has:

\begin{eqnarray*}
\bfGamma^*&=& D_s\otimes I\,\bfGamma\\[1mm]
&=&[D_s\hat\P_s^\top\hat\Omega\otimes I]\, S( \hat\bfuno\otimes y_0 + h\hat\I_s\otimes I\,\bfGamma)\, [\hat\P_s \P_s^\top\Omega\otimes I] \,\nabla H(\bfuno\otimes y_0 + h\I_s\otimes I\,\bfGamma)\\[1mm]
&=&[\hat\P_s^\top\hat\Omega P_{k_1}\otimes I]\, S( \hat\bfuno\otimes y_0 + h\hat\I_s\otimes I\,\bfGamma)\, [P_{k_1}^2\hat\P_s \P_s^\top\Omega\otimes I] \,\nabla H(\bfuno\otimes y_0 + h\I_s\otimes I\,\bfGamma)\\[1mm]
&=&[\hat\P_s^\top\hat\Omega \otimes I]\, S( P_{k_1}\hat\bfuno\otimes y_0 + hP_{k_1}\hat\I_s D_s^2\otimes I\,\bfGamma)\, [P_{k_1}\hat\P_s D_s^2\P_s^\top\Omega\otimes I] \,\nabla H(\bfuno\otimes y_0 + h\I_s\otimes I\,\bfGamma)\\[1mm]
&=&[\hat\P_s^\top\hat\Omega \otimes I]\, S( \hat\bfuno\otimes y_0 + h(\hat\bfuno\bfe_1^\top -\hat\I_s) \otimes I\,\bfGamma^*)\, [\hat\P_s \P_s^\top\Omega P_{k_2}\otimes I] \,\nabla H(\bfuno\otimes y_0 + h\I_s\otimes I\,\bfGamma)\\[1mm]
&=&[\hat\P_s^\top\hat\Omega \otimes I]\, S( \hat\bfuno\otimes (y_0+h\hat\Gamma_0) -h\hat\I_s \otimes I\,\bfGamma^*)\, [\hat\P_s \P_s^\top\Omega\otimes I] \,\nabla H(P_{k_2}\bfuno\otimes y_0 + hP_{k_2}\I_sD_s^2\otimes I\,\bfGamma)\\[1mm]
&=&[\hat\P_s^\top\hat\Omega \otimes I]\, S( \hat\bfuno\otimes y_1 -h\hat\I_s \otimes I\,\bfGamma^*)\, [\hat\P_s \P_s^\top\Omega \otimes I] \,\nabla H(\bfuno\otimes y_0 + h(\bfuno\bfe_1^\top-\I_s)\otimes I\,\bfGamma^*)\\[1mm]
&=&[\hat\P_s^\top\hat\Omega \otimes I]\, S( \hat\bfuno\otimes y_1 -h\hat\I_s \otimes I\,\bfGamma^*)\, [\hat\P_s \P_s^\top\Omega \otimes I] \,\nabla H(\bfuno\otimes (y_0 + h\hat\Gamma_0)-h\I_s\otimes I\,\bfGamma^*)\\[1mm]
&=&[\hat\P_s^\top\hat\Omega \otimes I]\, S( \hat\bfuno\otimes y_1 -h\hat\I_s \otimes I\,\bfGamma^*)\, [\hat\P_s \P_s^\top\Omega \otimes I] \,\nabla H(\bfuno\otimes y_1-h\I_s\otimes I\,\bfGamma^*).
\end{eqnarray*}
Consequently, the statement follows.\,\QED

\subsection{Solving the discrete problem}\label{blend}
Preliminarily, let us recall that, by the properties of Legendre polynomials (see, e.g., \cite{LIMbook2016}), one has:
\begin{equation}\label{Xs}
\hat\P_s^\top\hat\Omega\hat\I_s = \P_s^\top\Omega\I_s = X_s \equiv \pmatrix{cccc} \xi_0 &-\xi_1\\ \xi_1 &0 &\ddots\\ &\ddots &\ddots &-\xi_{s-1}\\&&\xi_{s-1}&0\endpmatrix, \quad \xi_i = \frac{1}{2\sqrt{|4i^2-1|}},
\end{equation}
and
\begin{equation}\label{POP1}
\hat\P_s^\top\hat\Omega\hat\P_s = I_s, \qquad \P_s^\top\Omega\bfuno = \bfe_1.
\end{equation}
Next, we observe that the discrete problem (\ref{bfGampro}) naturally induces the fixed-point iteration:
\begin{equation}\label{fixit}
\bfGamma^{\ell+1}=[\hat\P_s^\top\hat\Omega\otimes I]\, S( \hat\bfuno\otimes y_0 + h\hat\I_s\otimes I\,\bfGamma^\ell)\, [\hat\P_s \P_s^\top\Omega\otimes I] \,\nabla H(\bfuno\otimes y_0 + h\I_s\otimes I\,\bfGamma^\ell),\qquad \ell=0,1,\dots.\quad
\end{equation}
The following theorem provides sufficient condition for its convergence.

\begin{theo}\label{fixconv} Assume that the function $S(y)$ and $\nabla H(y)$ at the right-hand side in (\ref{poisson}) are Lipschitz with constant $L$ in a suitable closed ball of radius $\rho$ centered at $y_0$. Let us denote
$$\alpha_0 = \max_y\{ \|S(y)\|: \|y-y_0\|\le\rho\}, \quad \alpha_1 = \max_y\{ \|\nabla H(y)\|: \|y-y_0\|\le\rho\},\quad \alpha = \max\{\alpha_0,\alpha_1\},$$ where, hereafter, $\|\cdot\|$ denotes the infinity norm. Moreover, let us set:
$$\rho_1 = \|\hat\P_s^\top\hat\Omega\|\cdot\|\hat\P_s\|\cdot \|\P_s^\top\Omega\|, \qquad \rho_2 = \|\hat\I_s\|+\|\I_s\|.$$
Then, starting at $\bfGamma^0=\bfzero$, the iteration (\ref{fixit}) converges, provided that the timestep $h$ satisfies
\begin{equation}\label{hsmall}
0\le h<h^*:=\min\left\{\frac{1}{\rho_1\rho_2\alpha L},~\frac{\rho}{\rho_2\|S(y_0)\nabla H(y_0)\|}\right\}.
\end{equation}
\end{theo}
\proof Starting from $\bfGamma^0=\bfzero$, by virtue of (\ref{POP1}) one obtains $\bfGamma^1=\bfe_1\otimes S(y_0)\nabla H(y_0)$.
Consequently, from (\ref{hsmall}) it follows that
$$\|h\hat\I_s\otimes I\bfGamma^1\|,\,\|h\I_s\otimes I\bfGamma^1\|<\rho.$$
Next, assume that
$$\|h\hat\I_s\otimes I\bfGamma^i\|,\,\|h\I_s\otimes I\bfGamma^i\|<\rho, \qquad i=\ell-1,\ell.$$
Then,
\begin{eqnarray*}
\|\bfGamma^{\ell+1}-\bfGamma^\ell\| &=& \|[\hat\P_s^\top\hat\Omega\otimes I]\, S( \hat\bfuno\otimes y_0 + h\hat\I_s\otimes I\,\bfGamma^\ell)\, [\hat\P_s \P_s^\top\Omega\otimes I] \,\nabla H(\bfuno\otimes y_0 + h\I_s\otimes I\,\bfGamma^\ell)\\[1mm]
&&-[\hat\P_s^\top\hat\Omega\otimes I]\, S( \hat\bfuno\otimes y_0 + h\hat\I_s\otimes I\,\bfGamma^{\ell-1})\, [\hat\P_s \P_s^\top\Omega\otimes I] \,\nabla H(\bfuno\otimes y_0 + h\I_s\otimes I\,\bfGamma^{\ell-1})\|\\
&=&\|[\hat\P_s^\top\hat\Omega\otimes I]\, S( \hat\bfuno\otimes y_0 + h\hat\I_s\otimes I\,\bfGamma^\ell)\, [\hat\P_s \P_s^\top\Omega\otimes I] \,\nabla H(\bfuno\otimes y_0 + h\I_s\otimes I\,\bfGamma^\ell)\\[1mm]
&&\pm [\hat\P_s^\top\hat\Omega\otimes I]\, S( \hat\bfuno\otimes y_0 + h\hat\I_s\otimes I\,\bfGamma^\ell)\, [\hat\P_s \P_s^\top\Omega\otimes I] \,\nabla H(\bfuno\otimes y_0 + h\I_s\otimes I\,\bfGamma^{\ell-1})\\[1mm]
&&-[\hat\P_s^\top\hat\Omega\otimes I]\, S( \hat\bfuno\otimes y_0 + h\hat\I_s\otimes I\,\bfGamma^{\ell-1})\, [\hat\P_s \P_s^\top\Omega\otimes I] \,\nabla H(\bfuno\otimes y_0 + h\I_s\otimes I\,\bfGamma^{\ell-1})\|\\[1mm]
&\le& h\rho_1\rho_2\alpha L \|\bfGamma^\ell-\bfGamma^{\ell-1}\|.
\end{eqnarray*}
Consequently, we have a contraction, for all $h<h^*.\,\QED$\medskip

As is clear from (\ref{hsmall}), when the right-hand side in (\ref{poisson}) has a large linear part, the fixed-point iteration (\ref{fixit}) may require a quite small timestep to converge. When this happens,  a Newton-type iteration for solving (\ref{bfGampro}) would be more appropriate, which we now sketch below.

To begin with, we observe that, when $k_1=k_2=s$, then $\hat\P_s=\P_s$, $\hat\I_s=\I_s$,  $\hat\Omega=\Omega$,  which are $s\times s$ matrices, and $\P_s^{-1}=\P_s^\top\Omega$ (see (\ref{POP1})). In such a case, (\ref{bfGampro}) becomes (using the notation in  (\ref{poisson}))
$$G(\bfGamma) := \bfGamma - \P_s^\top\Omega\otimes I f(\bfuno\otimes y_0+h\I_s\otimes I\,\bfGamma) = \bfzero.$$
The use of the simplified Newton iteration for solving such an equation then becomes, by virtue of (\ref{Xs}):
\begin{eqnarray}\nonumber
\mbox{solve:}&& [I_s\otimes I-h X_s\otimes f'(y_0)]\Delta^\ell = -G(\bfGamma^\ell)\\[1mm]
\mbox{set:}&&\bfGamma^{\ell+1} = \bfGamma^\ell+\Delta^\ell, \qquad \ell=0,1,\dots,\label{newt}
\end{eqnarray}
where, as is usual, $f'(y_0)$ denotes the Jacobian of $f(y)$ evaluated at $y_0$. Also in this case, the initial guess $\bfGamma^0=\bfzero$ can be conveniently used. However, since the dimension of the discrete problem is always $s$, whichever are $k_1$ and $k_2$, we shall continue using the iteration (\ref{newt}) even when $k_1>s$ and/or $k_2>s$, so that $G(\bfGamma^\ell)$ is now evaluated according to (\ref{bfGampro}). In any case, the straight solution of the linear systems in (\ref{newt}) would require the factorization of a matrix whose dimension is $s$ times larger than that of the continuous problem (\ref{poisson}). In order to obtain a comparably effective, though less expensive, Newton-type iteration, we shall consider a corresponding Newton-splitting {\em blended iteration}. This iteration, at first devised in \cite{B2000}, has been studied in \cite{BM2002} (see also \cite{BM2009}), and implemented in the computational codes {\tt BiM} \cite{BM2004} and {\tt BiMD} \cite{BMM2006}, respectively solving stiff ODE-IVPs and linearly implicit DAEs. Later on, it has been considered for HBVMs \cite{BIT2011,BFCI2014,BMR2018}, as well as for other line integral methods (see, e.g., \cite{BMR2019}). We here only sketch the final iteration:\footnote{The parameter $\rho_s$ is determined according to a linear analysis of convergence \cite{BM2002,BM2009} and, as usual, $\sigma(X_s)$ denotes the spectrum of matrix $X_s$.}
\begin{eqnarray}\nonumber
\mbox{set:}&&\rho_s = \min_{\lambda\in\sigma(X_s)}|\lambda|, \quad \Theta = [I-h\rho_s f'(y_0)]^{-1}, \quad \bfGamma^0=\bfzero\\[1mm] \label{blend}
\mbox{for $\ell=0,1,\dots$:}&&\bfeta^\ell = -G(\bfGamma^\ell), \quad \bfeta_1^\ell = [\rho_sX_s^{-1}]\otimes I\,\bfeta^\ell\\[1mm]\nonumber
&&\bfGamma^{\ell+1} = \bfGamma^\ell - I_s\otimes\Theta\left[\bfeta_1^\ell + I_s\otimes(I-\Theta)(\bfeta^\ell-\bfeta_1^\ell)\right]
\end{eqnarray}
Consequently, the iteration (\ref{blend}) now requires only the factorization of one matrix (i.e., $\Theta^{-1}$), having the same size as that of the continuous problem.

\section{Numerical tests}\label{num}
We here report a few numerical tests describing particular instances of the gyrocenter dynamics of a charged particle \cite{ZHTZ2015}, i.e. :
\begin{enumerate}
\item the case of a dipole magnetic field {  with zero electric potential},

\item the case of a tokamak magnetic field  {  with zero electric potential},

\item {  the case of a dipole magnetic field with quadratic electric field.}
\end{enumerate}

{  In the first two problems, the electric potential $\phi(x)$ is assumed to be 0, which is called the {\em temporal gauge}.} Consequently, the magnetic vector potential (\ref{A}) is enough to completely define the problem, along with the constant $\mu$ and the initial conditions. All  numerical tests have been done on an Intel i7 computer with 16GB of memory, running Matlab 2019a. 

In the first two problems, the fixed-point iteration (\ref{fixit}) is used for solving the generated discrete problems, since the linear part of the corresponding right-hand sides turns out to be very small (of the order of unity),  so that there is no gain in using the blended iteration (\ref{blend}), whose computational cost per iterate is higher. {  However, in the third test problem the blended iteration turns out to have a superior performance.} Moreover, in the numerical tests we set $k_1=s$ and $k_2\equiv k\ge s$, so that all LIM$(s,k_2,s)$ methods have order $2s$.

\begin{table}[t]
\caption{Maximum Hamiltonian error when solving the dipole magnetic field problem on the interval $[0,10^3]$, by using LIM$(s,k,s)$ with timestep $h=0.4$.} \label{tab1}
\centerline{
\begin{tabular}{|r|rrrrr|}
\hline
      &\multicolumn{5}{|c|}{$s$}\\
      \hline
 $k$ & 1 & 2 & 3 & 4 & 5\\
 \hline
   1 & 2.689e-02 &  &  &  &  \\ 
  2 & 6.163e-04 & 5.103e-03 &  &  &  \\ 
  3 & 3.549e-06 & 5.551e-05 & 2.785e-04 &  &  \\ 
  4 & 8.366e-08 & 6.909e-07 & 8.613e-06 & 1.374e-05 &  \\ 
  5 & 1.425e-09 & 1.371e-08 & 1.040e-07 & 3.796e-07 & 6.394e-07 \\ 
  6 & 3.256e-11 & 4.590e-10 & 1.998e-09 & 7.869e-09 & 1.552e-08 \\ 
  7 & 1.776e-15 & 8.698e-12 & 5.307e-11 &  1.455e-10  & 2.828e-10 \\ 
  8 & 2.220e-15 & 1.776e-15 & 5.653e-13 &  2.850e-12 & 4.602e-12 \\ 
  9 &                  & 1.776e-15 & 1.776e-15 &  1.776e-15 &  1.776e-15 \\ 
  10&                 &                  &  2.220e-15 & 1.776e-15  & 1.776e-15 \\ 
 \hline
\end{tabular}}
\end{table}
\begin{table}[t]
\caption{Dipole magnetic field problem solved on the interval [0,40] with timestep $h$ (times in {\em sec}).}\label{tab2}
\scriptsize
\centerline{
\begin{tabular}{|r|rrr|rrr|rrr|rrr|rrr|rrr|}
\hline
                 &\multicolumn{3}{|c|}{LIM(1,7,1)} &\multicolumn{3}{|c|}{LIM(2,8,2)}  
                 &\multicolumn{3}{|c|}{LIM(3,9,3)} &\multicolumn{3}{|c|}{LIM(4,9,4)} 
                 &\multicolumn{3}{|c|}{LIM(5,9,5)} \\
       \hline
$h$ & err & rate & time & err & rate & time & err & rate & time & err & rate & time & err & rate & time \\
\hline
$0.4$
&1.05e\,00 &--- &  0.2 &1.58e-02 &--- &  0.3 &1.82e-03 &--- &  0.2 &4.12e-05 &--- &  0.2 &1.78e-07 &--- &  0.2 \\ 
$2^{-1}0.4$
&2.90e-01 &1.9 &  0.3 &1.71e-03 &3.2 &  0.4 &3.35e-05 &5.8 &  0.4 &9.44e-08 &8.8 &  0.3 &1.68e-09 &6.7 &  0.3 \\ 
$2^{-2}0.4$
&7.44e-02 &2.0 &  0.4 &1.20e-04 &3.8 &  0.6 &5.16e-07 &6.0 &  0.6 &4.74e-10 &7.6 &  0.5 &1.11e-12 &10.6 &  0.5 \\ 
$2^{-3}0.4$
&1.87e-02 &2.0 &  0.6 &7.69e-06 &4.0 &  0.9 &8.06e-09 &6.0 &  1.0 &2.18e-12 &7.8 &  0.9 &4.00e-13 &** &  0.9 \\ 
$2^{-4}0.4$
&4.68e-03 &2.0 &  1.0 &4.84e-07 &4.0 &  1.6 &1.26e-10 &6.0 &  1.6 &3.32e-13 &** &  1.6 &                &    &      \\ 
$2^{-5}0.4$
&1.17e-03 &2.0 &  1.8 &3.03e-08 &4.0 &  2.8 &1.54e-12 &6.4 &  3.0 &                &         &           &                &         &           \\ 
$2^{-6}0.4$
&2.93e-04 &2.0 &  3.1 &1.89e-09 &4.0 &  5.0 &1.34e-13 &** &  5.4 &                &         &           &                &         &           \\ 
$2^{-7}0.4$
&7.32e-05 &2.0 &  5.8 &1.19e-10 &4.0 &  9.2 &                &         &           &                &         &           &                &         &           \\ 
$2^{-8}0.4$
&1.83e-05 &2.0 & 10.5 &7.02e-12 &4.1 & 17.0 &                &         &           &                &         &           &                &         &           \\ 
$2^{-9}0.4$
&4.58e-06 &2.0 & 20.0 &  7.91e-13      &  **     & 32.9 &                &       &         &                &       &         &                &       &         \\ 
$2^{-10}0.4$
&1.14e-06 &2.0 & 37.5 &                &       &         &                &       &         &                &       &         &                &       &         \\ 
$2^{-11}0.4$
&2.86e-07 &2.0 & 68.2 &                &       &         &                &       &         &                &       &         &                &       &         \\ 
\hline
\end{tabular}}
\end{table}

\begin{figure}[t]
\centerline{\includegraphics[height=6.5cm,width=8cm]{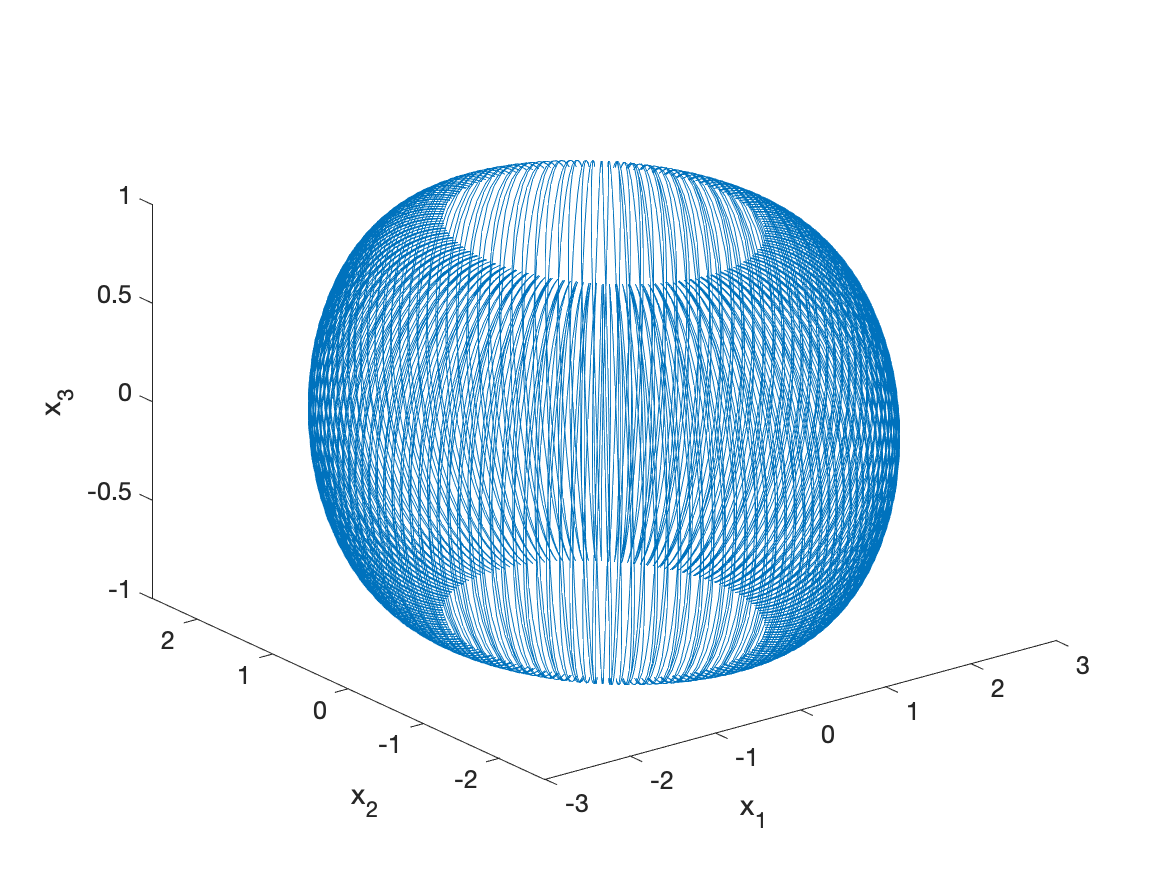}\qquad
\includegraphics[height=6cm,width=7cm]{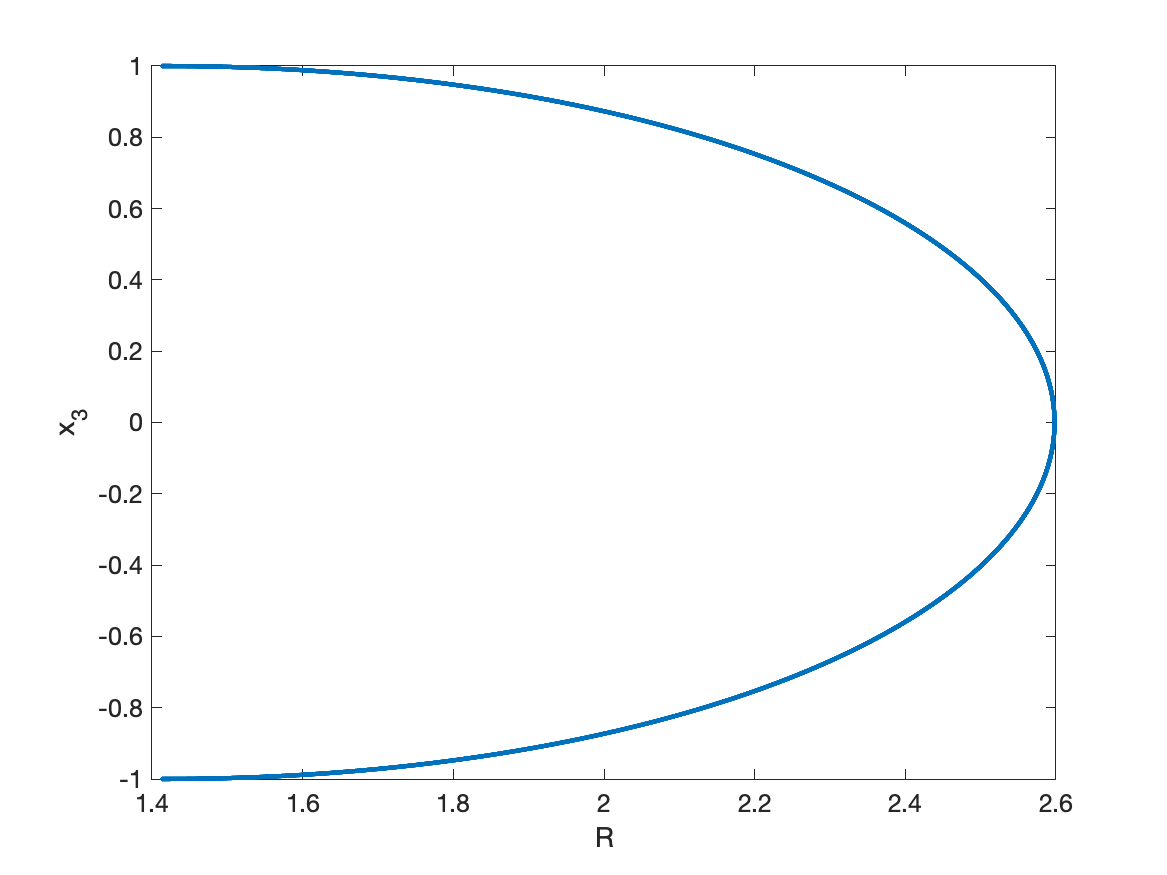}}
\caption{Dipole magnetic problem. Left plot: solution trajectory in the phase space; right plot: solution in the coordinates (\ref{rteta}).}\label{pic1}
\end{figure}

\paragraph{Dipole magnetic field.}\label{dipo}
The dipole magnetic field has wide applications in different branches of physics, underlying many examples of
cosmic magnetic fields, such as the earth magnetic field and the neutron stars magnetic field.
The vector potential is 
\begin{equation}\label{Ad}
A(x) = \frac{M}{\rho^3}\pmatrix{ccc} x_2, &  -x_1, & 0\endpmatrix^\top, \qquad \rho=\|x\|,
\end{equation}
with $M$ being the dipole moment, which can be either positive or
negative.\footnote{In the case of earth, the constant $M=-8\times10^{15}$.} Consequently, one obtains:
\begin{eqnarray}\nonumber
B(x) &=& -\frac{M}{\rho^5}\pmatrix{ccc} 3x_1x_3, & 3x_2x_3, & 2x_3^2-x_1^2-x_2^2\endpmatrix^\top,\\[1mm] \label{Bd}
\|B(x)\| &=& |M|\frac{\sqrt{\rho^2+3x_3^2}}{\rho^4}, \\[1mm] \nonumber
b(x) &=& \frac{-M/|M|}{\rho\sqrt{\rho^2+3x_3^2}}\pmatrix{ccc} 3x_1x_3, & 3x_2x_3, & 2x_3^2-x_1^2-x_2^2\endpmatrix^\top.
\end{eqnarray}
In the numerical tests, we set
\begin{equation}\label{Mmud}
M=10^3, \qquad \mu=10^{-2}, \qquad y(0) = (\,1,\,1,\,1,\,0.01\,)^\top.
\end{equation}
In Figure\,\ref{pic1} we show the 3D trajectory of the gyrocenter  (left plot), and its 2D representation (right plot)  in the coordinates
\begin{equation}\label{rteta}
R =\sqrt{x_1^2+x_2^2}, \qquad x_3,
\end{equation}
for $t\in[0,10^3]$.

At first, in Table~\ref{tab1} we show the maximum Hamiltonian error, by using a timestep $h=0.4$, for the LIM$(s,k,s)$ methods, $s=1,\dots,5$, and $k=s,\dots,10$. As one may expect, as the value of $k$ increases, the Hamiltonian error decreases. In particular, in order to obtain energy-conservation up to round-off error level, a value $k=7$ is enough when $s=1$, $k=8$ is sufficient in the case $s=2$, whereas $k=9$ is appropriately chosen, for $s=3,4,5$.

Then, in Table~\ref{tab2}, we list the solution error, along with the estimate convergence rate, for the LIM(1,7,1), LIM(2,8,2), and LIM$(s,9,s)$, $s=3,4,5$, methods when using a timestep $h=2^{-i}0.4$, $i\ge0$, to cover the interval $[0,40]$. We also list the corresponding execution times (in {\em sec}). As one may see, the convergence order $2s$ is satisfied for all methods. Moreover,  the higher-order methods turns out to be much more efficient than the lower-order ones. We notice that the LIM$(1,7,1)$, which is the lowest-order method, is equivalent to the method used in \cite{ZLQT2019}.

\begin{figure}[t]
\centerline{\includegraphics[height=6cm,width=8cm]{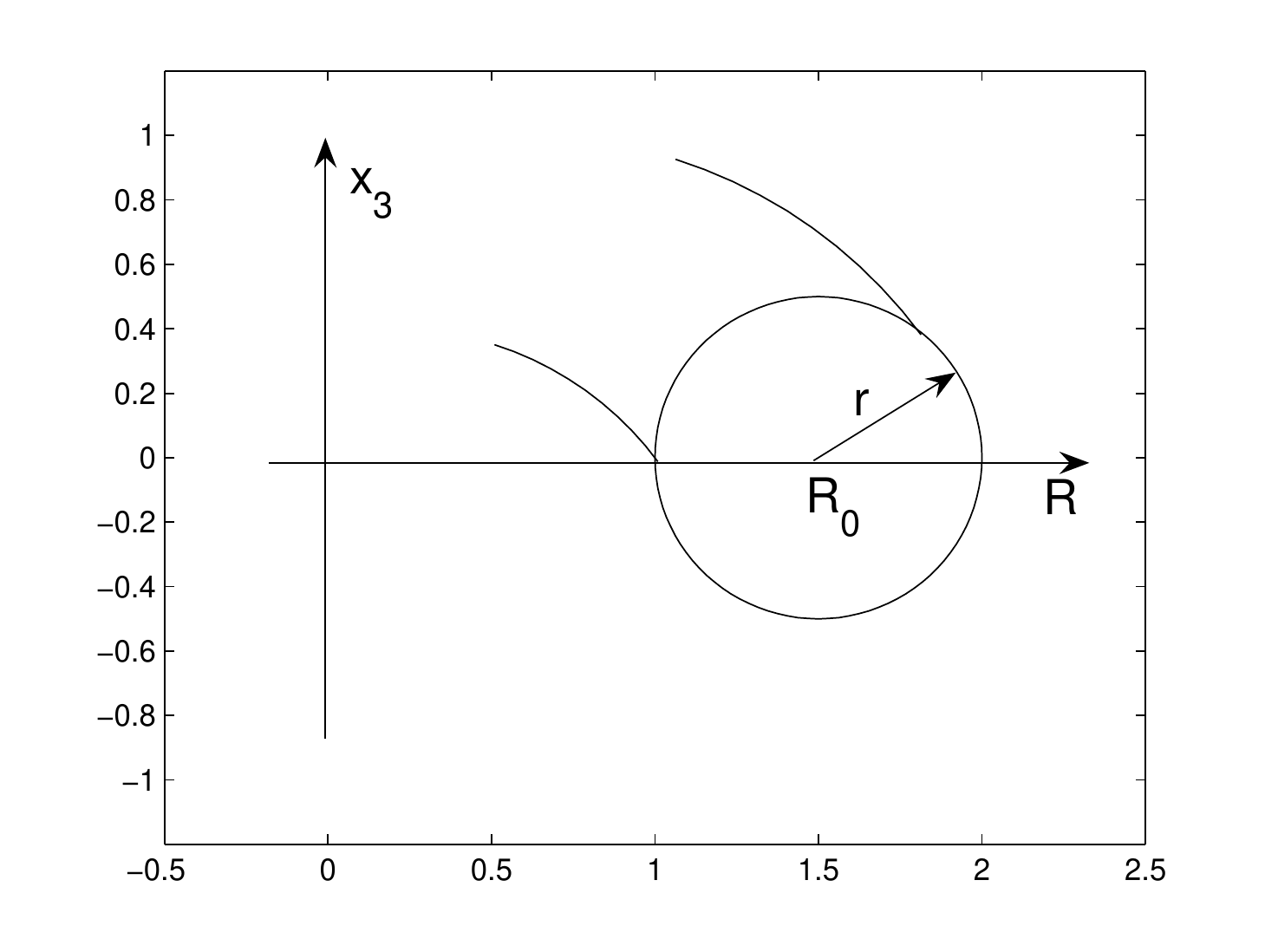}}\caption{2D tokamak geometry with circular concentric flux surfaces \label{pic2}}
\end{figure}

\paragraph{Tokamak magnetic field.}\label{toka} 

\begin{figure}[p]
\centerline{\includegraphics[height=6.5cm,width=8cm]{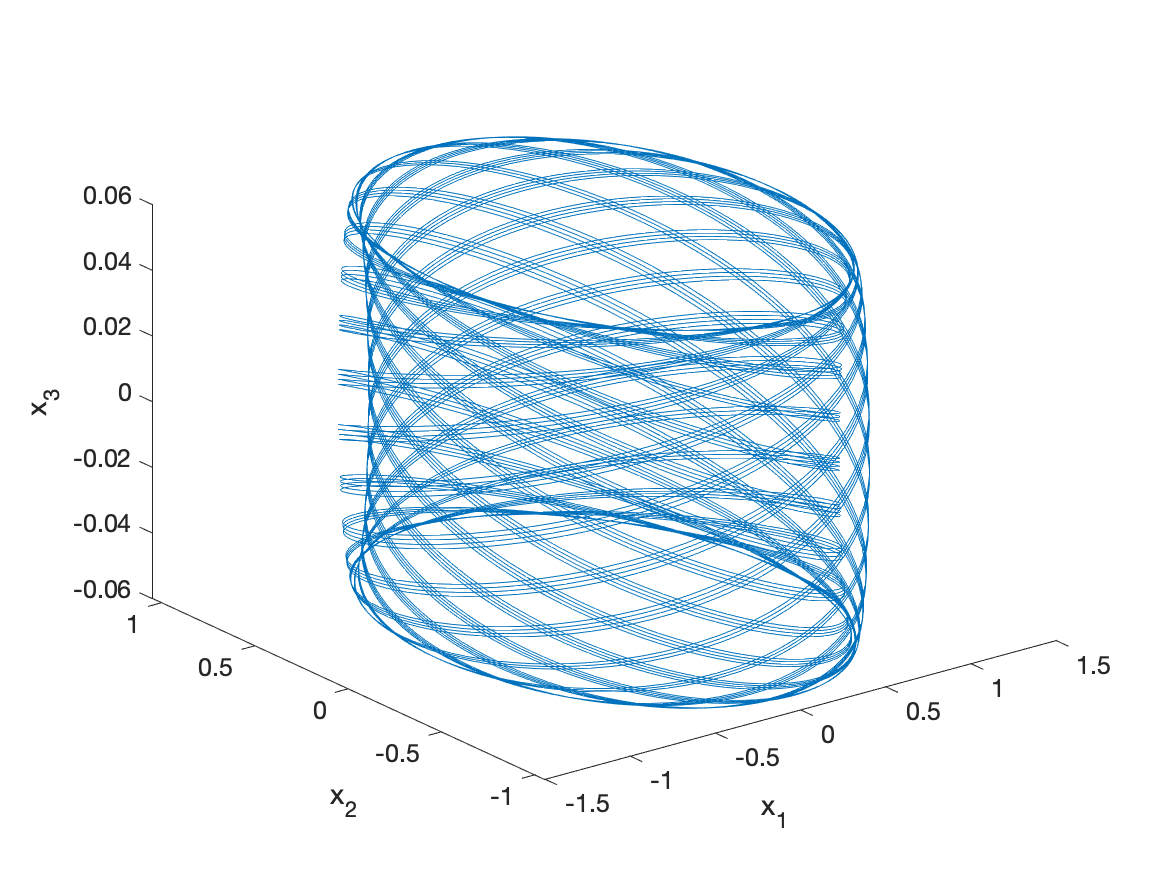}
\includegraphics[height=6cm,width=7cm]{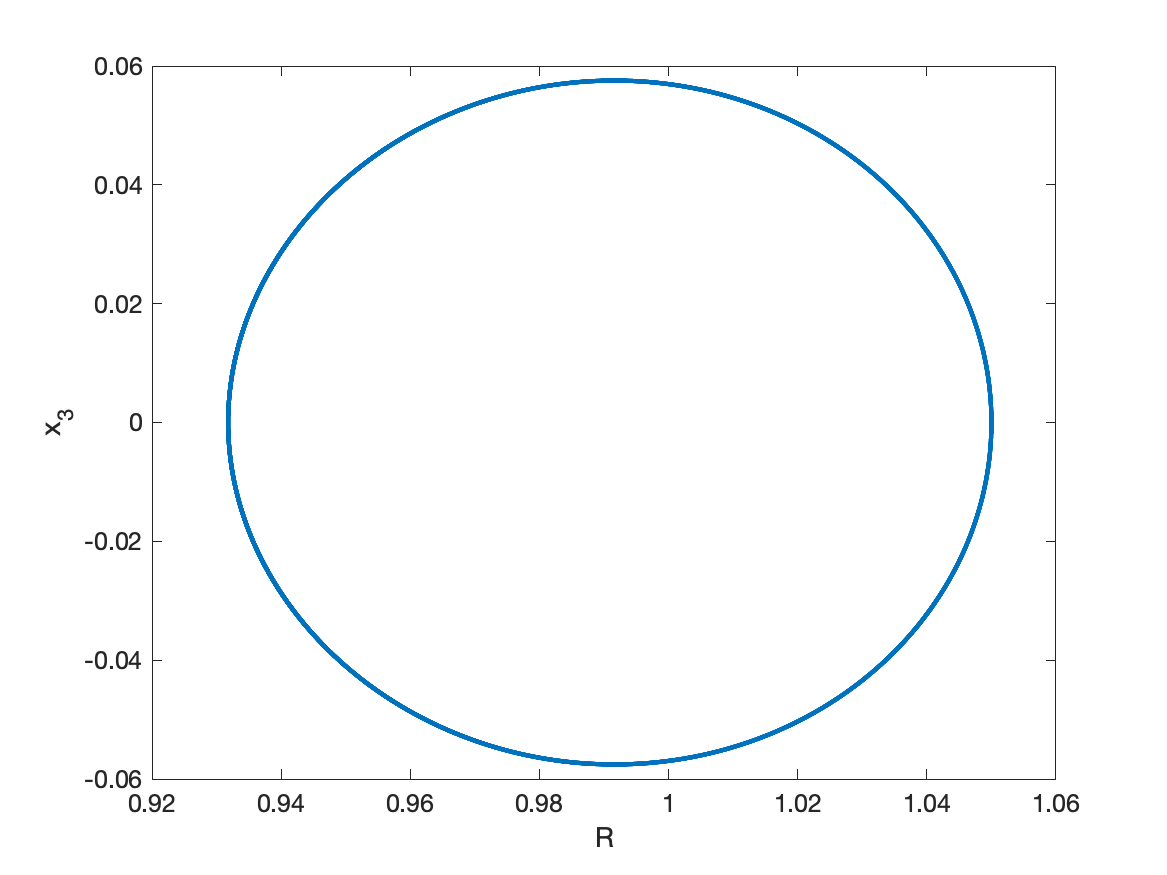}}
\caption{Transit orbit. Left plot: solution trajectory in the phase space; right plot: solution in the coordinates (\ref{rteta}).}\label{pic3}

\bigskip
\bigskip
\centerline{\includegraphics[height=6.5cm,width=8cm]{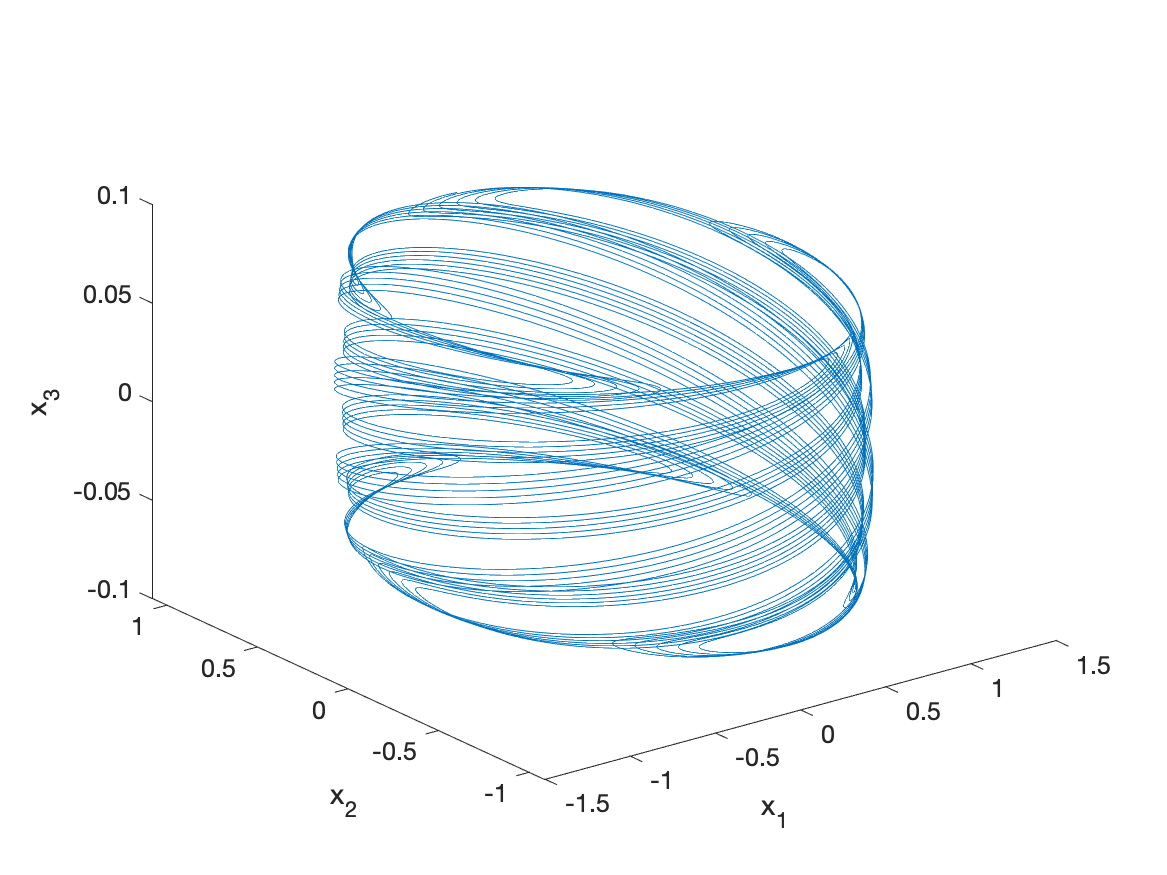}\qquad
\includegraphics[height=6cm,width=7cm]{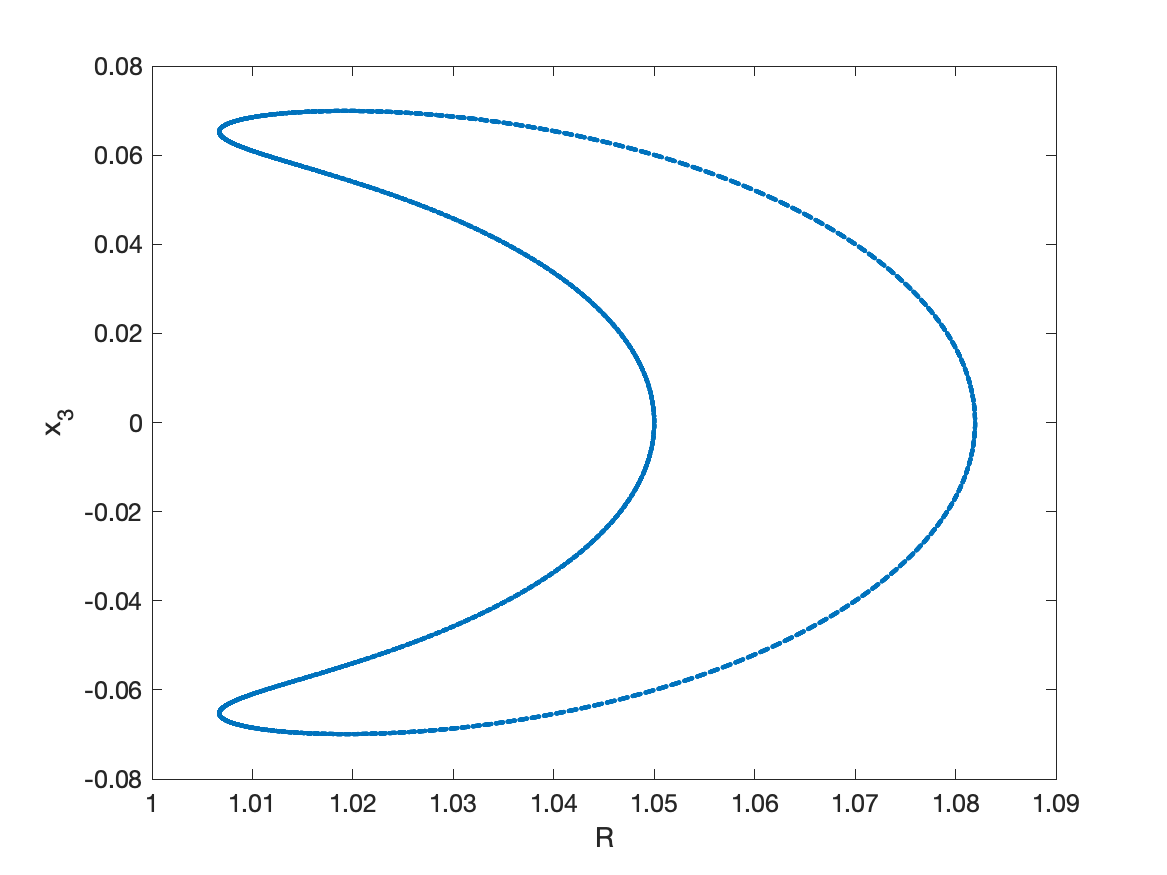}}
\caption{Banana orbit. Left plot: solution trajectory in the phase space; right plot: solution in the coordinates (\ref{rteta}).}\label{pic4}
\end{figure}

The tokamak magnetic field is used to confine a hot plasma in the shape of a torus. The axisymmetric tokamak geometry is illustrated in Figure\,\ref{pic2}, where $R$ is defined according to (\ref{rteta}), $r = \sqrt{(R-R_0)^2+x_3^2}$, and $R_0$ is the main radius.
We give an example for the vector potential 
\begin{equation}\label{At}
A(x) = \frac{B_0}{2q R^2}\pmatrix{ccc} qR_0x_1x_3-x_2r^2, & qR_0x_2x_3+x_1r^2, &
-qR^2R_0\log(R/R_0)\endpmatrix^\top,
\end{equation}
where $B_0$ is the magnetic field on the axis and the constant $q$ is
the safety factor, from which one obtains
\begin{eqnarray}\nonumber
B(x) &=& \frac{B_0}{qR^2}\pmatrix{ccc}
-x_1x_3-qR_0x_2, &
-x_2x_3+qR_0x_1, &
R(R-R_0)\endpmatrix^\top,\\[1mm] \label{Bt}
\|B(x)\| &=& \frac{B_0}{qR}\sqrt{r^2+q^2R_0^2},\\[1mm] \nonumber
b(x) &=&\frac{1}{R\sqrt{r^2+q^2R_0^2}}\pmatrix{ccc}
-x_1x_3-qR_0x_2, &
-x_2x_3+qR_0 x_1, &
R(R-R_0)\endpmatrix^\top,
\end{eqnarray}
In the numerical tests, we shall choose the parameters as follows:
\begin{equation}\label{tokapar}
R_0 = 1, \qquad B_0 = 1, \qquad q=2, \qquad \mu = 2.25\cdot 10^{-6}.
\end{equation}
Moreover, we shall choose the initial condition as (see (\ref{y}))
\begin{equation}\label{y0}
y_0 = \pmatrix{cccc} 1.05, &0, &0, & u_0\endpmatrix^\top,
\end{equation}
by considering the following two values of $u_0$:
\begin{itemize}
\item $u_0 = 0.0008117$, generating a transit orbit, i.e.,  a circular orbit, in the coordinates (\ref{rteta}).

\item $u_0 = 0.0004306$, which generates a banana shaped orbit, in the coordinates (\ref{rteta});
\end{itemize}
The 3D trajectories in the interval $[0,10^6]$ are respectively shown in the left plots of Figures\,\ref{pic3}--\ref{pic4}, whereas those in the coordinates (\ref{rteta}) are shown in the corresponding right plots.

We now further show that the higher-order LIMs are much more efficient than the lower-order ones. In fact,  by selecting $k$ large enough so that energy-conservation is granted even for relatively large stepsizes, e.g. $k=20$, the methods can be regarded as {\em spectral methods in time} for large values of $s$ (but smaller than $k$). Spectral methods in time proved to be very effective in many instances, such as highly-oscillatory problems \cite{BMR2018},  stiff-oscillatory problems \cite{BIMR2019}, fractional equations \cite{ABI2019}, and other problems \cite{BMR2019}.\footnote{The reasons for their effectiveness has been studied in \cite{ABI2018}.} For this purpose, let us solve the transit orbit problem on the interval $[0,10^8]$, by using a timestep as large as $h=8\cdot10^3$, and the banana orbit problem, on the same interval, by using a timestep $h=10^4$. In Table~\ref{tab3} we list the obtained results by using the LIM$(s,20,s)$ method, $s=1,\dots,16$. For both problems, a reference solution has been computed by using the LIM(18,20,18) method with the same timestep. In the table, one has, for each selected value of $s$:
\begin{itemize}
\item *** if the nonlinear iteration (\ref{fixit}) does not converge; 
\item the total number of nonlinear iterations (\ref{fixit}) for covering the iteration interval;
\item the corresponding execution time (in {\em sec});
\item the maximum error w.r.t. the reference solution. 
\end{itemize}
As one may see, only for $s$ large enough the problems are solved. Moreover, even though each iteration (\ref{fixit}) has a cost which increases with $s$, nevertheless, the higher $s$, the smaller the number of iterations needed, so that the overall execution time decreases with $s$, and so does the error. It is worth noticing that this effectiveness is made possible provided that a vector function is used for evaluating $S(y)$ and $\nabla H(y)$, as in the present case. In fact, this allows to exploit the vector architecture of nowadays processors.\footnote{A similar remark was observed in \cite[Remark\,3]{BMR2018}.} The obtained results clearly testify the effectiveness of higher-order energy-conserving LIMs w.r.t. lower-order ones. 

\begin{table}
\caption{Numerical results when solving the transit orbit ($h=8\cdot 10^3$) and the banana orbit ($h=10^4$) tokamak problems on the interval $[0,10^8]$ by using the energy-conserving LIM$(s,20,s)$ methods (*** means non convergence of the nonlinear iteration, times are in {\em sec}).}
\label{tab3}
\centerline{\begin{tabular}{|r|rrr|rrr|}
\hline
       &\multicolumn{3}{c|}{transit orbit}&\multicolumn{3}{c|}{banana orbit}\\
\hline
$s$ & iterations & time & error & iterations & time & error\\
\hline
  1 &        *** &  *** & *** &        *** &  *** & *** \\ 
  2 &        *** &  *** & *** &        *** &  *** & *** \\ 
  3 &        *** &  *** & *** &        *** &  *** & *** \\ 
  4 &        *** &  *** & *** &        *** &  *** & *** \\ 
  5 &        *** &  *** & *** &        *** &  *** & *** \\ 
  6 &        *** &  *** & *** &        *** &  *** & *** \\ 
  7 &        *** &  *** & *** &        *** &  *** & *** \\ 
  8 &        *** &  *** & *** &   773705 & 84.7 & 3.1e\,00 \\ 
  9 &  1018824 & 114.4 & 3.0e\,00 &   570191 & 64.4 & 8.5e-01 \\ 
 10 &   734527 & 85.0 & 1.2e\,00 &   494422 & 57.2 & 6.2e-02 \\ 
 11 &   625527 & 74.2 & 1.1e-01 &   457523 & 54.2 & 1.6e-02 \\ 
 12 &   569554 & 68.6 & 9.2e-03 &   436163 & 52.6 & 1.3e-03 \\ 
 13 &   533843 & 65.7 & 7.1e-04 &   419205 & 51.7 & 1.9e-04 \\ 
 14 &   509484 & 63.9 & 5.0e-05 &   410197 & 51.7 & 1.1e-05 \\ 
 15 &   501218 & 64.1 & 2.5e-06 &   402775 & 51.4 & 1.3e-06 \\ 
 16 &   493683 & 63.9 & 8.1e-07 &   399053 & 51.6 & 2.5e-07 \\ 
\hline
\end{tabular}}
\end{table}

\paragraph{Dipole magnetic field with quadratic electric field.}\label{dipo1}
We now consider a further problem, aimed at showing the potentialities of the blended iteration (\ref{blend}). This problem is still defined by (\ref{Ad})--(\ref{Mmud}), with the initial condition replaced by
\begin{equation}\label{y0fi} y(0) = (\,1,\,1,\,0.01,\,0.01\,)^\top. \end{equation} Moreover, now the electric potential defining the Hamiltonian (\ref{H}) is not zero, and is given by
\begin{equation}\label{fi}
\phi(x) = \frac{1}2 x^\top G x, \qquad G = \diag\left(\,1,\,1,\,10^4\,\right). 
\end{equation}
According to the result of Theorem~\ref{fixconv}, now the fixed-point iteration may encounter stepsize restrictions. In Table~\ref{tab4} we list the obtained results by solving the problem on the smallest interval containing $[0,10^3]$, commensurable with the (approximately) maximum  timestep $h_{max}$ allowed by the used iteration, either fixed-point (f-p) or blended (blend), using the same methods considered in Table~\ref{tab2}. In such a case, as is clear, we are not discussing accuracy but, instead, we want to emphasize the robustness of the nonlinear iteration.\footnote{Remarkably enough, all methods remain still energy-conserving, even when using large timesteps.} 

\begin{table}[t]
\caption{Dipole magnetic field with quadratic electric field problem solved on an interval containing $[0,10^3]$; (approximate) maximum timestep ($h_{max}$) allowed by the iteration (either fixed-point (f-p) or blended (blend)), total number of iterations (it), mean number of iterations per step (it/step), and execution times (in {\em sec}).}\label{tab4}

\medskip
\centerline{
\begin{tabular}{|r|rr|rr|rr|rr|rr|rr|}
\hline
                 &\multicolumn{2}{|c|}{LIM(1,7,1)} &\multicolumn{2}{|c|}{LIM(2,8,2)}  
                 &\multicolumn{2}{|c|}{LIM(3,9,3)} &\multicolumn{2}{|c|}{LIM(4,9,4)} 
                 &\multicolumn{2}{|c|}{LIM(5,9,5)} \\
                  \hline
iteration    &   f-p & blend &   f-p & blend &   f-p & blend &   f-p & blend &   f-p & blend \\
\hline
$h_{max}$ &   .01        &   47     & .02           & 72    &  .04          & 86     & .05          & 103     & .06 & 120\\
it                & 4329357 &  880    &  2697596 & 1120& 5013330  & 1333 & 3225342 &  1420  &2713788 & 1599\\ 
it/step        & 43.3        &  40      &  54.0        & 80    & 200.5       & 111.1 & 161.3     &  142.0 & 162.8     &  177.6\\
time           & 282.4      &    0.1   &  307.3      & 0.1   &  605.7      & 0.2    & 392.0      & 0.2      &  345.0& 0.2\\
\hline
\end{tabular}}
\end{table}

From the results listed in Table~\ref{tab4}, one may observe that, for LIM$(s,k,s)$ methods used for solving this problem:
\begin{itemize}
\item the maximum allowed timestep, $h_{max}$, increases with $s$;
\item the mean number of iterations per step increases with $h_{max}$ (and, then, with $s$); 
\item the blended iteration allows using much larger (indeed, huge) timesteps, w.r.t. the fixed-point iteration, resulting in much smaller execution times.
\end{itemize}

\subsection{Conclusions}\label{fine}
In this paper we have developed arbitrarily high-order methods for Poisson problems, with a major emphasis on the  simulation of the gyrocenter dynamics of a charged particle in a constant magnetic field, which is a relevant problem in plasma physics. The methods are derived and studied within the framework of {\em line integral methods}, and their efficient implementation has been also sketched. The reported numerical tests duly confirms the theoretical achievements. As a future direction of investigations, we plan to study the speed-up of the convergence of the nonlinear iteration, when large value of $s$ are taken into account, corresponding to the use of the methods as  spectral methods in time.

\subsection*{Acknowledgements} The research of the third author was supported by the National Key Research and Development Program
 (Grant No.\,2017YFE0301700).

\end{document}